\documentclass[twoside]{article}
\usepackage[accepted]{aistats2020}
\usepackage{hyperref}
\usepackage{epsfig}
\usepackage{xcolor}
\usepackage{amsthm}
\usepackage{amssymb}
\usepackage{amsfonts}
\usepackage{amsmath}
\usepackage{amsopn}
\usepackage{graphicx}
\usepackage{epstopdf}
\usepackage{algorithmic}
\usepackage{algorithm}
\usepackage{paralist}
\usepackage{tikz}
\usepackage{float}
\usepackage{wrapfig}
\usepackage{ulem}
\usepackage{stmaryrd}
\usepackage{booktabs}
\usepackage{colortbl}
\usepackage{appendix}
\usepackage{afterpage}
\usepackage{longtable}
\usepackage{varwidth}
\usepackage{caption}
\usepackage{listings}
\usepackage{multirow}
\usepackage{url}
\usepackage{svg}
\usepackage{pstricks}
\usepackage[noabbrev]{cleveref}

\usepackage{rbase}
\usepackage{macros}

\newcommand{\inner}[1]{\left< #1 \right>}
\theoremstyle{remark}
\newtheorem*{remark}{Remark}
\theoremstyle{definition}

\theoremstyle{plain}

\usepackage[round]{natbib}
\usepackage{amsmath,amssymb,mathtools}

\begin{document}



\twocolumn[

\aistatstitle{Asymmetric Multiresolution Matrix Factorization}

\aistatsauthor{ Pramod Kaushik Mudrakarta \And Shubhendu Trivedi \And  Risi Kondor }

\aistatsaddress{ University of Chicago \And  CSAIL, MIT \And Flatiron Institute, NY \\ University of Chicago} ]

\begin{abstract}
	Multiresolution Matrix Factorization (MMF) was recently introduced as an alternative to the dominant low-rank paradigm in order to capture structure in matrices at multiple different scales. Using ideas from multiresolution analysis (MRA), MMF teased out hierarchical structure in symmetric matrices by constructing a sequence of wavelet bases. While effective for such matrices, there is plenty of data that is more naturally represented as nonsymmetric matrices (e.g. directed graphs), but nevertheless has similar hierarchical structure. In this paper, we explore techniques for extending MMF to any square matrix. We validate our approach on numerous matrix compression tasks, demonstrating its efficacy compared to low-rank methods. Moreover, we also show that a combined low-rank and MMF approach, which amounts to removing a small global-scale component of the matrix and then extracting hierarchical structure from the residual, is even more effective than each of the two complementary methods for matrix compression. 

%With a growing number of data having hierarchical structure (e.g. social network graphs, circuit simulations, etc.), there is a strong need for developing methods that efficiently represent them for downstream tasks such as matrix compression. A critical challenge is that matrices that represent such data are often not low rank. Multiresolution Matrix Factorization was first developed in order to extend harmonic analysis to discrete spaces, and consequently model hierarchical structure by teasing out discrete wavelet transforms from symmetric matrices. While this method has been successfully used for a range of applications, there is plenty of data that cannot be represented by symmetric matrices (e.g. directed graphs). In recent literature, there is no consensus on how hierarchical structure manifests in nonsymmetric matrices. In this paper, we explore a few practical techniques for extending Multiresolution Matrix Factorization to any square matrix. Our experiments on matrix compression show that the proposed techniques can be more effective than low-rank-based methods on a variety of matrices.
\end{abstract}
\section{Introduction}

A wide variety of naturally occurring data, represented as large matrices, often has complex hierarchical structure. For instance, this is the case in data originating from social networks, complex biological interactions, circuit simulations and so on. With the increasing size of such matrices, there is a need to represent and compress them efficiently. Towards this end, a number of methods have been proposed~\citep{mahoney2009cur,coifman2006diffusion,hackbusch1999sparse, halko2009}. All such methods make an implicit or explicit assumption about the nature of such matrices, which serves as a prior to facilitate efficient factorization and compression. Frequently, the assumption in question is that the matrices are of low rank. For example, Principal Component Analysis (PCA) decomposes a symmetric matrix $A$ into the form \begin{align*}A = Q^T D\ts Q\end{align*} where $Q$ is a orthogonal matrix and $D$ is diagonal. Notably, such methods are single-level algorithms. Moreover, the columns of $Q$ are usually dense, which does not adequately capture strong locality properties that might be exhibited, for instance, in real world graphs. 

Taking a different tack, Multiresolution Matrix Factorization (MMF)~\citep{kondor2014multiresolution} framed matrix factorization as a type of multiresolution analysis (MRA) over a discrete metric space whose metric structure is represented by a self-adjoint operator. Thus, in contrast to low-rank approaches, MMF factorizes a symmetric matrix $A$ in the following multilevel manner: \begin{align*}A = Q_1^T Q_2^T \dots Q_L^T H\ts Q_L \dots Q_2 Q_1\end{align*} where each of the carefully selected $Q_i$ matrices are orthogonal and extremely sparse, and the matrix $H$ is close to diagonal. MMF is able to discover soft hierarchical arrangements of localized structures in the data represented by the matrix.  It is amenable to parallelization and fast implementations~\citep{kondor2015parallel}. The wavelet bases uncovered by MMF also provide a natural sparse approximation for functions defined over the space. For instance, if the matrix $A$ is the normalized graph Laplacian, then the factorization provides us with wavelet bases for approximating functions on graphs.

While MMF has been shown to be useful in various practical settings~\citep{teneva2016multiresolution,ding2017multiresolution,ithapu2017decoding}, one notable drawback of the method is that it is wholly defined using a self-adjoint operator such as the graph Laplacian. Therefore, MMF can only handle a multi-level factorization of symmetric matrices. However, many naturally occurring data are represented more concisely as an nonsymmetric matrix; the underlying space might only have a quasi-metric structure. This is the case for instance, when we are working with directed graphs. For data living on such spaces, hierarchical or multiscale structure is as important, albeit more complicated, and thus there is much to be gained by defining appropriate sparse wavelet bases for function approximation.

Defining such wavelet bases, in particular for functions on directed graphs, is an active area of research~\citep{sevi2018harmonic, satoshi2019, mhaskar2018}. However, there is currently no consensus on the correct approach for doing so. Further, existing methods are difficult to implement in practice for downstream applications. In this paper, we do not make an attempt to directly extend multiresolution analysis (MRA) in the spirit of MMF to such spaces, but rather provide practical methods to leverage intuitions used in MMF to make it operable on general square matrices. We show experimentally that our methods can perform very well in various real world matrix compression tasks. 

To conclude this section, our contributions in this paper are twofold: 1) we provide two practical ways to adapt MMF to nonsymmetric matrices, 2) we argue that extending MRA to nonsymmetric matrices is not analogous to developing a multi-level SVD; in particular, we show experimentally that it is important to consider interactions between different levels in the hierarchical structure and lastly 3) we also show that a combined low-rank and MMF approach, which amounts to removing a small global-scale component of the matrix and then extracting hierarchical structure from the residual, is even more effective than each of the two complementary methods for matrix compression

After a description of notation used in this paper, in~\Cref{sec:mmf}, we give a brief introduction to Multiresolution Matrix Factorization.

\section{Notation}
We denote $[n] = {1, 2, \ldots, n}$. Inner products between vectors are denoted by the operator $\inner{\cdot,\cdot}$. The $(i,j)$-th entry of a matrix $A$ is denoted by $A_{i,j}$. If $S\subseteq [n]$ is a set of indices, $A_{S,S}$ represents the submatrix of $A$ corresponding to those indices. $\overline{S}$ denotes the complement of a set. When $S$ can be interpreted as an ordered set, for instance, if it is a set of indices, then $S[i]$ represents the $i$-th index from that set. 
\section{Multiresolution Matrix Factorization}\label{section:mmf}
\label{sec:mmf}

Given a measurable space $X$ and square integrable functions $L_2(X)$ on $X$, multiresolution analysis (MRA) constructs a sequence of spaces of functions of increasing smoothness:
\[
\def\arraystretch{1.5}
\begin{array}{ccccccccc}
L_2(X) = V_0 & \rightarrow & V_1 &\rightarrow &V_2 &\rightarrow &\ldots &\rightarrow &V_L \\
& \searrow & & \searrow & & \searrow & & \searrow & \\
& & W_1 & & W_2 & & \ldots & & W_L 
\end{array}\] where each $V_\ell$ is split into a smoother space $V_{\ell+1}$ and a rougher space $W_{\ell+1}$, such that $V_\ell = V_{\ell+1} \oplus W_{\ell+1}$. The length-scale over which a function $f\in V_\ell$ varies increases with $\ell$ and thus makes this decomposition an analysis of $L_2(X)$ multiscale/multiresolution.

\cite{kondor2014multiresolution} extends MRA to functions on discrete spaces (e.g. undirected graphs) by interpreting a given symmetric matrix as a discrete version of a continuous self-adjoint operator. Specifically, for $\vert X \vert = n$,~\cite{kondor2014multiresolution} defines multiresolution analysis with respect to a symmetric matrix $A\in\R^{n\times n}$ as a successive filtration \begin{align}\label{eq:filtration}V_L \subset \ldots \subset V_1 \subset V_1 \subset V_0 = L(X) \cong \R^n\end{align} where each $V_\ell$ has an orthonormal basis $\Phi_\ell := \lbrace \phi_m^\ell \rbrace$, and its complementary space $W_\ell$ has an orthonormal basis $\Psi_\ell := \lbrace \psi_m^\ell \rbrace$ that satisfy the following three axioms:
\begin{description}
    \item[MRA1] The sequence~\cref{eq:filtration} is a filtration of $\R^n$ in terms of smoothness. In particular, \[ \eta_\ell = \sup_{v\in V_\ell} \frac{\inner{v,Av}}{\inner{v,v}} \] decays with $\ell$.
    \item[MRA2] The basis vectors of $W_\ell$ are localized in the sense that \[ \mu_l = \max_{m\in\lbrace 1,\ldots,d_l\rbrace} \norm{\psi_m^\ell}_0 \] increases no faster than a certain rate.
    \item[MRA3] Suppose the matrix $U_\ell$ represents the matrix expressing $\Psi_l \cup \Phi_l$ in terms of the previous basis $\Phi_{\ell-1}$. Then, $U_\ell$ is sparse. ($\Phi_0$ is taken to be the standard basis.)
\end{description}

Once an overcomplete multiscale basis satisfying MRA1-3 for $\R^n$, i.e., $U_\ell, \ell=1,\ldots, L$, is determined, it can be used to compress the matrix $A$. Applying each of the $U_\ell$'s sequentially, we obtain versions of $A$ compressed to various scales: 

\[A\mapsto \underbrace{U_1 A\ts U_1^T}_{A_1} \mapsto \ldots 
\mapsto \underbrace{U_L\ldots U_2 U_1 A\ts U_1^T U_2^T\ldots U_L^T}_{A_L = H}\]   

Since the choice of the $U_\ell$'s aligns with MRA1,~\cite{kondor2014multiresolution} argues that the matrices $A_1,\ldots, A_L$ and $H$ are approximately in core-diagonal form. That is, up to a permutation, each of those matrices has a block-diagonal form, with one block being a diagonal matrix. 

Formally, Multiresolution Matrix Factorization is defined as an approximate multilevel factorization of the form 
\begin{align}
\label{eq:mmf}
A \approx Q_1^T Q_2^T \ldots Q_L^T H\tts Q_L \ldots Q_2 Q_1, 
\end{align}
where the matrices $Q_1,\ldots,Q_L$ and $H$ obey the following conditions:
\begin{enumerate}[~1.]
    \item Each $Q_\ell \in \R^{n\times n}$ is orthogonal and highly sparse. 
    In the simplest case, each $Q_\ell$ is a Givens rotation, i.e., a matrix which differs from the identity in 
    just the four matrix elements
    \begin{align*}
    [Q_\ell]_{i,i} = \cos \theta,&\qqquad [Q_\ell]_{i,j} = - \sin \theta,\\
    [Q_\ell]_{j,i} = \sin \theta,&\qqquad [Q_\ell]_{j,j} = \cos \theta,
    \end{align*}
    for some pair of indices $(i,j)$ and rotation angle $\theta$. 
    Multiplying a vector with such a matrix rotates it counter-clockwise by $\theta$ in the $(i,j)$ plane. 
    More generally, $Q_\ell$ is a so-called 
    $k$-point rotation, 
    which rotates not just two, but $k$ coordinates. 
    \item Typically, in MMF factorizations $L\<=O(n)$, and the size of the active part of the \m{Q_\ell} matrices decreases according to a set schedule $n = \delta_0 \ge \delta_1 \ge \ldots \ge \delta_L$. 
    More precisely, there is a nested sequence of sets $[n] = S_0 \supseteq S_1 \supseteq \ldots \supseteq S_L$ 
    such that the $[Q_\ell]_{\overline{S_{\ell-1}},\overline{S_{\ell-1}}}$ 
    part of each rotation is the $n-\delta_{\ell-1}$ 
    dimensional identity. $S_\ell$ is called the \textbf{active set} at level $\ell$. In the simplest case,  $\delta_\ell=n-\ell$.
    \item $H$ is an $S_L$-core-diagonal matrix, which means that it   
    is block diagonal with two blocks: $H_{S_L,S_L}$, called the core, which is dense, and 
    $H_{\wbar{S_L},\wbar{S_L}}$ which is diagonal. In other words, $H_{i,j}\!=\!0$ unless $i,j\tin S_L$ or $i\!=\!j$.  
\end{enumerate}

Given a symmetric matrix, the factorization is carried out by minimizing the Frobenius norm error between the original matrix and its factorized form. Algorithms such as \textsc{GreedyJacobi}~\citep{kondor2014multiresolution} or pMMF~\citep{teneva2016multiresolution} are typically used in practical applications~\citep{teneva2016multiresolution,ding2017multiresolution,mudrakarta2017generic}

Having discussed the core ideas behind MMF, we briefly review some related art before moving forward.

%%%%%%%%%%%%%%%%%%%%%%%%%%%%%%%%%%%%%%%%%%%%%%%%%%%%%%%%%%%%%%%%%%%%%%
%%%%%%%%%%%%%%%%%%%%%%%%%%%%%%%%%%%%%%%%%%%%%%%%%%%%%%%%%%%%%%%%%%%%%%
%%%%%%%%%%%%%%%%%%%%%%%%%%%%%%%%%%%%%%%%%%%%%%%%%%%%%%%%%%%%%%%%%%%%%%
%%%%%%%%%%%%%%%%%%%%%%%%%%%%%%%%%%%%%%%%%%%%%%%%%%%%%%%%%%%%%%%%%%%%%%
%%%%%%%%%%%%%%%%%%%%%%%%%%%%%%%%%%%%%%%%%%%%%%%%%%%%%%%%%%%%%%%%%%%%%%

\ignore{
\section{Multiscale structure in nonsymmetric matrices}

In this section, I want to write about why its not trivial to extend MMF to nonsymmetric matrices, and also explain a little bit about the intuition behind the formulations in the next sections. Two points to address:
\begin{itemize}
    \item How well is MRA1 satisfied in the formulations proposed? i.e., what measure of smoothness are we using? 
    \item 
\end{itemize}

\section*{Notes (just for guidance - not intended to be part of the paper)}
\begin{itemize}
    \item The broad goal here is to efficiently compress and estimate large matrices arising in the real world
    \item This process involves making assumptions on the matrix -- the most common one being the low-rank assumption, which is typically modeled via PCA
    \item The problem with PCA is that the eigenvectors are always dense, which does not capture cases where the matrix ``more closely couples certain clusters of nearby coordinates than those farther apart w.r.t. some underlying topology''. Sparse PCA solves the locality problem, but does not capture global scale structure.
    \item MMF is a method that aims to capture multiscale structure in discrete spaces such as graphs.  Matrices are considered as operators on vectors, and MMF essentially decomposes that operator in terms of a multiscale basis. 
    \item Given a continuous operator such as a Laplacian. The application of this operator on a function can be expressed in terms of the eigenfunctions of the operator. In the case of the Laplacian, the eigenfunctions are the Fourier basis functions. 
    
    Similarily, when the symmetric matrix is the operator, the application of this operator (i.e., matrix-vector product) can be expressed in terms of eigenvectors of the matrix. Thus, in this way the eigenvalue decomposition is an analogue of the Fourier decomposition. 
    \item Multiresolution Analysis (MRA) captures multiscale structure in a given space of functions by repeatedly reducing the space so that the smoothness of functions increases in some sense. For functions on $\mathbb{R}$, Mallat (1989) defines MRA in terms of the wavelets constructed via dilation and translation of certain canonical functions. 
    
    \pramod{Question: how does Mallat show that this form of constructing wavelets correspond to the idea of compressing the space by smoothness of functions? In particular, what exactly is the measure of smoothness that he uses, and how is it related to symmetricity/asymmetricity of an operator that is used to determine the smoothness of a function?}
    \item The MMF paper defines MRA on discrete spaces similarly, in terms of successively compressing a given sucspace. Note that each subpsace here would be defined by a finite set of basis vectors. Three axioms are defined that charactersize the smoothness of each space (MRA1), the localization (MRA2), and the existence of fast wavelet transform (MRA3). 
    
    \pramod{What does it mean to relax the criterion that each $V_\ell$ must be a strict subset of $V_{\ell-1}$?}
    
    \pramod{How is smoothness defined in the asymmetric case?}
    
    \item Now that we can determine a wavelet basis in matrix form $U_1, \ldots U_L$ , the given operator $A$ can be transformed into it by performing 
    $$H = U_L\ldots U_1 A\ts U_1^T \ldots U_L^T$$
    
    In the symmetric case, the wavelet functions are assumed to be close to eigenvectors of $A$, so that the MRA1 criterion ensures that each subspace is smoother than the previous. This makes the structure of $H$ core-diagonal. MMF-factorizable matrices are defined to be those for which applying the wavelet transformation results in $H$ being core-diagonal. 
    
    \pramod{In reality, $H$ never is core-diagonal. We manually zero out entries of the wavelet-transformed matrix so that it comes to that form. That makes MMF an approximate factorization.}
    
    \pramod{In some of the asymmetric versions we develop, this criterion is relaxed, thereby violating the MRA1 axiom a bit, but that results in better compression.}
\end{itemize}
}
\section{Related Work}

The work most closely related to some of the ideas in this paper is that on defining wavelets on functions on directed graphs. Two notable recent works on the theme are that of \cite{sevi2018harmonic} and \cite{mhaskar2018}. \cite{sevi2018harmonic} did so by using a random walk operator, while \cite{mhaskar2018} used an extension of diffusion polynomial frames~\citep{mhaskarmaggioni}, using a different operator. In fact, the influential work on diffusion wavelets~\citep{coifman2006diffusion}, provides a more general theory in which symmetry of the diffusion operator is not necessary. However, it is unclear how to lift these ideas for defining a principled asymmetric MMF that has also a clear MRA interpretation. On the side of matrix factorization, works that have some resemblance to our approach are \cite{LiYangButterfly, LiYangInterpolate}, in which the left and right matrices don’t have to be related by transposition. However these approaches don't have a clear wavelet based motivation.
\section{MMF on nonsymmetric matrices}

Based on the MMF approach described in~\Cref{section:mmf}, in the sequel we provide two approaches that extend MMF to general square matrices. 

\subsection{MMF via an additive decomposition}

We begin by observing that any square matrix $A\in\R^{n\times n}$ can be written as the sum of a symmetric matrix and a skew-symmetric matrix as follows
\begin{align}
\label{eq:square_decomp}
A = \underbrace{\frac{A + A^T}{2}}_\text{symmetric} + \underbrace{\frac{A - A^T}{2}}_\text{skew-symmetric} 
\end{align}

We propose performing MRA on each term separately. Since the first term is symmetric, it directly yields to a multiresolution matrix factorization. In particular, the factorization of $\frac{A + A^T}{2}$ can be obtained via a parallel variant of MMF such as pMMF~\citep{kondor2015parallel}. Thus, the problem of computing a multiresolution factorization of an asymmetric square matrix is now reduced to computing a multiresolution factorization for a skew-symmetric matrix. 

\subsubsection{MMF for skew-symmetric matrices}
For skew-symmetric matrices, the MRA axioms, especially MRA1, are no longer applicable. The reason being that under the definition of smoothness in MRA1 $\sup_{v\in V\subseteq\R^n} \frac{\inner{v,Av}}{\inner{v,v}}$, for a skew-symmetric matrix $A$, the smoothness is always zero.

In symmetric MMF, the choice of the sparse orthogonal factors and the diagonal part of the core-diagonal matrix $H$ are based on the assumption that the multiscale basis vectors determined by the MMF algorithm are as close to eigenvectors as possible. In other words, the filtration defined in~\cref{eq:filtration} is a sequence of subspaces that is as close to a spectral filtration as possible. Indeed, MMF approximately diagonalizes a symmetric matrix using $O(n)$ rotations, where $n$ is the dimension of the matrix. We consider this intuition as our starting off point to extend MMF to skew-symmetric matrices.

Note that a matrix $A\in\R^{n\times n}$ is skew-symmetric if $A^T = -A$. Skew-symmetric matrices have complex eigenvalues and are therefore not diagonalizable in real space. Thus, a skew-symmetric matrix can never be perfectly represented using the same form of MMF as for a symmetric matrix over reals. In order to overcome this issue, we consider a well known observation~\citep{murnaghan1931canonical} that a skew-symmetric matrix, under a unitary transformation, assumes the so-called Murnaghan canonical form,\begin{align}
\begin{bmatrix}
\Lambda_1 & 0 & \ldots & 0\\
0 & \Lambda_2 & \ldots & 0\\
\vdots & \vdots & \ddots & \vdots  \\
0 & \ldots & 0 & \Lambda_L
\end{bmatrix}
\label{eqn:murnaghan}
\end{align}
where $\Lambda_i, i=1,\ldots, L$ are $2\times 2$ matrices of the form \[\Lambda_i = \begin{pmatrix} 0 & \lambda_i \\ -\lambda_i & 0\end{pmatrix}\]
where $\lambda_i \in \R$. 

This observation affords a multiresolution factorization for skew-symmetric matrices identical to that for symmetric matrices with the difference that the core-diagonal matrix is allowed to have the following form:
\begin{align}
H = \begin{bmatrix}
H_\text{core} & 0 & \ldots & \ldots & 0\\
0 & \Lambda_1 & 0 & \ldots & 0\\
0 & 0 & \Lambda_2 & \ldots & 0\\
\vdots & \vdots & \ddots & \ddots & \vdots \\
0 & \ldots & \ldots & 0 & \Lambda_L
\end{bmatrix}
\label{eqn:asymmetric_core_diagonal}
\end{align}
where $\Lambda_i, i=1,\ldots, L$ are $2\times 2$ matrices as defined above.

A couple of observations based on this form are in order. First,  notice that the skew-symmetricity of the matrix in its MMF is always preserved.  Second, this form also allows for interaction between wavelets at adjacent levels. This feature is distinct from the original MMF procedure in which no such interactions occur. To obtain the factorization, the procedure is similar to the \textsc{GreedyJacobi} algorithm from~\cite{kondor2014multiresolution}, with the difference that in in the last step, instead of zeroing out off-diagonal elements, we sparsify the matrix into the Murnaghan normal form as described above.

Thus, the multiresolution decomposition of a square matrix can be done by factorizing each of the two parts from~\cref{eq:square_decomp} separately and writing $A$ as a sum of two MMFs. In other words, we may write 
\begin{align*}
A \approx& P_1^T\ldots P_L^T H_\text{sym}\ts P_L \ldots P_1  \\ & +  Q_1^T\ldots Q_L^T H_\text{skew}\ts Q_L \ldots Q_1
\end{align*} 
where 
\begin{align}
\label{eqn:asymmetric_MMF}
\frac{A + A^T}{2} &\approx P_1^T\ldots P_L^T H_\text{sym} P_L \ldots P_1 \\
\frac{A - A^T}{2} &\approx Q_1^T\ldots Q_L^T H_\text{skew} Q_L \ldots Q_1
\end{align}

where $H_\text{skew}$ has the form defined in~\cref{eqn:asymmetric_core_diagonal}.

In the next section, we present another approach for multiresolution factorization.

\subsection{MMF via direct factorization}

In this approach, we aim to capture hierarchical structure by operating on the square matrix directly rather than decomposing it into an additive form. Recall that we would like to identify clusters of coordinates that a matrix $A$ more closely couples than others, and determine an appropriate multiresolution analysis of $\R^{n}$. When the matrix is symmetric, the linear transformation that it represents is self-adjoint. In the asymmetric case, depending on whether we are transforming $x\rightarrow x^TA$ or $x \rightarrow Ax$, we may end up with different ways in which $A$ might couple the coordinates. 

We consider two different ways of performing multiresolution analysis of $\R^n$, one that is based on $A$, and another on $A^T$. Suppose $P\in\R^{n\times n}$ represents the wavelet transform in the former case, and $Q\in\R^{n\times n}$ of the latter. Thus, we have 

\[ H = P^TAQ\]

which is the matrix $A$ in the bases $P$ and $Q$. In the symmetric case, $H$ takes a core-diagonal form, based on the fact that the eigenvalue decomposition converts a symmetric matrix to diagonal, and MMF uses orthonormal transformations that are designed to be as close to eigenvectors as possible. In the asymmetric case, we relax this assumption and allow some of the off-diagonal elements of $H$ to be nonzero. This is motivated by the fact that the frequencies or smoothness coefficients filtered by the left and right hand sides may not align with each other.

\begin{remark}
    When $A$ represents the adjacency matrix of a directed graph, we are essentially identifying hierarchical cluster structures considering 1) only the outgoing edges, and 2) only the incoming edges. The matrix $H$ gives a mapping between these two cluster structures.  
\end{remark}

We define asymmetric MMF to be a multiresolution factorization of the following form
\begin{align}
A = P_1P_2\ldots,P_L H Q_L^TQ_{L-1}^T\ldots Q_1^T
\end{align}, where
\begin{itemize}
    \item Each $P_i, Q_i, i=1,\ldots,L$ are highly sparse, orthogonal rotation matrices similar to the rotations in a symmetric MMF
    \item $H\in\R^{n\times n}$ is a core-sparse matrix, that is, under an appropriate permutation of row and columns, $H$ would be a $2\times 2$ block-diagonal matrix with the upper left block being dense, and the other 3 blocks would be highly sparse
    \item When $P_i = Q_i$ and $H$ is in core-diagonal form with a symmetric core, the above factorization would be a symmetric MMF. 
\end{itemize}

The computation of the asymmetric MMF is obtained using an approach similar to the symmetric MMF. The algorithm we propose is designed to minimize the Frobenius norm error 
\[\| A - P_1P_2\ldots,P_L H Q_L^TQ_{L-1}^T\ldots Q_1^T\|_{F} \]
We determine the left and right rotations in a greedy manner. That is, we compute a sequence of rotated version of $A$ 
\[A \rightarrow P_1^TAQ_1 \rightarrow \ldots \rightarrow P_L^T\ldots P_1^TAQ_1\ldots Q_L.\]
Each of the rotations is computed by identifying sets of similar rows (for the left rotations) and columns (for the right rotations), and determining an appropriate $k$-point rotation. This requires us to compute Gram matrices on both the rows as well as the columns. Finally, we sparsify the final rotated version of $A$ to determine the core-sparse matrix $H$. We explore three ways of sparsification
\begin{enumerate}
    \item \textsc{CoreDiagonal}: retain only the diagonal elements of $H$
    \item \textsc{TopN}: retain the top $n$ elements in terms of magnitude
    \item \textsc{GreedyTopN}: select $n$ elements by magnitude greedily such that no two of the selected elements fall in the same row/column
\end{enumerate}
The full procedure is summarized in~\Cref{alg:asymmetric_mmf}.

The running time of this algorithm is $O(n^3)$, similar to $\textsc{GreedyJacobi}$ routine used in the original symmetric MMF. Nevertheless, computation can be made faster by similarly extending the parallel MMF algorithm to maintain row and column Gram matrices and computing separate left and right rotations. It must be noted that applying the left rotations to $A$ does not affect the column Gram matrix, and applying the right rotations does not affect the row Gram matrix.

\begin{algorithm}
    \caption{Asymmetric MMF via direct factorization}
    \begin{algorithmic}
        \REQUIRE core size $d$, matrix $A_0=A\in\R^{n\times n}$
        \STATE Set of active rows $S^r_0 \gets [n]$
        \STATE Set of active columns $S^c_0 \gets [n]$
        \FOR{$\ell = 1$ \TO $n-d$}
        \STATE Select a random row $i$ from the set $S^r_\ell$
        \STATE $G \gets A_{S^r_\ell,S^c_\ell}A_{i,S^c_\ell}^T$ 
        \STATE $j \gets S^r_\ell[\arg\max(G)]$
        \STATE $P_\ell \gets \text{Jacobi rotation on rows } i,j$
        \STATE $A \gets P_\ell^T A$
        \STATE Remove $\arg\min_{t\in\lbrace i,j\rbrace}\norm{A_{t,S^c_\ell}}_2$ from $S^r_\ell$
        \STATE 
        \STATE Select a random column $i^\prime$ from the set $S^c_\ell$
        \STATE $G^\prime \gets A_{S^r_\ell,S^c_\ell}^TA_{S^r_\ell,i^\prime}$ 
        \STATE $j^\prime \gets S^c_\ell[\arg\max(G^\prime)]$
        \STATE $Q_\ell \gets \text{Jacobi rotation on rows } i^\prime,j^\prime$
        \STATE $A \gets Q_\ell$
        \STATE Remove $\arg\min_{t\in\lbrace i^\prime,j^\prime\rbrace}\norm{A_{S^r_\ell,t}}_2$ from $S^c_\ell$
        \ENDFOR
        \STATE $H = \textsc{CoreDiagonal}(A)$ or $\textsc{TopN}(A)$ or $\textsc{GreedyTopN}(A)$
        \RETURN $P_1,\ldots,P_{n-d}, H, Q_1, \ldots, Q_{n-d}$
    \end{algorithmic}
    \label{alg:asymmetric_mmf}
\end{algorithm}

\subsection{Numerical experiments}
\label{sec:asymmetric_experiments}
Similar to the evaluation of MMF in~\cite{teneva2016multiresolution}, we compare our asymmetric MMF formulations with the low-rank-based CUR decomposition~\citep{mahoney2009cur} for the task of compressing square matrices. The CUR decomposition is a widely used alternative to the singular value decomposition to compress large matrices. To reduce a matrix $A\in\R^{n\times n}$ to rank-$k$, the method proceeds by carefully selecting $k$ rows ($R\in\R^{n\times n}$) and columns ($C\in\R^{k\times n}$) and determines a matrix $U\in\R^{k\times k}$ such that the product $CUR$ is as close to the best rank-$k$ approximation of $A$ as possible. 

To set up a comparison, we sample 70 square matrices from the UFlorida sparse matrix repository~\citep{davis2011university}. The matrices in this repository were sourced from various fields of science and engineering and contain data arising from circuit simulations, social networks, chemical processes, 2D/3D problems, etc. All the sampled matrices were chosen such that the numerical symmetry was less than $25\%$, and had a maximum dimensionality of $100,000$.

\subsubsection{Results on MMF via an additive decomposition}

We first present results comparing the asymmetric MMF obtained by the additive decomposition to CUR. 
\Cref{tbl:mmf_cur_additive} shows the percentage of times the asymmetric MMF had a better compression error over CUR over 65 matrices. Since both the CUR and MMF algorithms are randomized, we report results averaged over 3 trials. CUR performs better than MMF in the majority of the cases. 

\begin{table*}[!htbp]
    \small 
    \centering 
    \begin{tabular}{lcccccc}
        \toprule 
        Matrix kind & num. matrices & 1\% & 10\% & 25\% & 50\% & 75\% \\
        \midrule
        \multicolumn{7}{c}{\textit{problems with underlying 2D/3D geometry}}\\
        structural problem&4&50.0&50.0&50.0&0.0&0.0\\
        computational fluid dynamics problem&3&0.0&0.0&0.0&0.0&0.0\\
        2D/3D problem&5&80.0&80.0&60.0&60.0&60.0\\
        materials problem&2&0.0&0.0&0.0&0.0&0.0\\
        \midrule 
        \multicolumn{7}{c}{\textit{problems with no underlying geometry}}\\
        economic problem&11&18.2&9.1&9.1&9.1&9.1\\
        directed graph&19&5.3&5.3&5.3&0.0&0.0\\
        power network problem &3&0.0&0.0&0.0&0.0&0.0\\
        circuit simulation problem&2&0.0&0.0&0.0&0.0&0.0\\
        chemical process simulation problem&23&0.0&0.0&8.7&0.0&0.0\\
        statistical/mathematical problem&1&100.0&0.0&0.0&0.0&0.0\\
        counter-example problem&1&0.0&0.0&0.0&0.0&0.0\\
        \midrule
        total wins: & &13.5  &10.8  &12.2  &5.4  &5.4 \\
        \bottomrule 
    \end{tabular}
    \caption{Matrix compression using asymmetric MMF (additive). Each entry of the rightmost 5 columns is the percentage of matrices on which MMF performed better than CUR, as measured by Frobenius norm compression error. Each column corresponds to the percentage size (w.r.t original size) to which a matrix is compressed.}
    \label{tbl:mmf_cur_additive}
\end{table*}

\subsubsection{Results on MMF via direct factorization}

Next, we report results on the direct factorization variant of the asymmetric MMF. In particular, we consider two different versions of the method based on what elements are retained in the core-sparse matrix $H$. In the first case, we retain the compressed $k\times k$ dense core, and $n-k$ elements chosen greedily according to magnitude and made sure that no two of these $n-k$ elements fall in the same row or column.

\Cref{tbl:mmf_cur_top_N_greedy} shows the percentage times this version of the asymmetric MMF wins over CUR in compression error. 
\begin{table*}[!htb]
    \small 
    \centering 
    \begin{tabular}{lcccccc}
        \toprule 
        Matrix kind & num. matrices & 1\% & 10\% & 25\% & 50\% & 75\% \\
        \midrule
        \multicolumn{7}{c}{\textit{problems with underlying 2D/3D geometry}}\\
        structural problem&4&50.0&50.0&25.0&0.0&0.0\\
        materials problem&1&0.0&0.0&0.0&0.0&0.0\\
        2D/3D problem&4&100.0&75.0&75.0&75.0&0.0\\
        computational fluid dynamics problem&3&100.0&100.0&0.0&0.0&0.0\\
        \midrule 
        \multicolumn{7}{c}{\textit{problems with no underlying geometry}}\\
        statistical/mathematical problem&1&100.0&100.0&100.0&0.0&0.0\\
        economic problem&9&77.8&22.2&11.1&33.3&11.1\\
        power network problem&2&100.0&0.0&50.0&0.0&0.0\\
        directed graph&15&80&20&6.7&0.0&0.0\\
        chemical process simulation problem&21&80.9&23.8&14.3&0.0&0.0\\
        circuit simulation problem&2&100.0&0.0&0.0&0.0&0.0\\
        counter-example problem&1&0.0&0.0&0.0&0.0&0.0\\
        \midrule
        total wins: & &79.4 &30.2 &17.5 &9.5 &1.6 \\
        \bottomrule 
    \end{tabular}
    \caption{Matrix compression using asymmetric MMF (\textsc{TopN}). See~\Cref{tbl:mmf_cur_additive} for a description of the table entries. }
    \label{tbl:mmf_cur_top_N_greedy}
\end{table*}
We observe that this version of MMF is significantly better than the additive version, and is able to capture structure in both problems with underlying known geometry and without underlying geometry. 

Relaxing the structure of the matrix $H$ further, we report the results for simply storing the top $n-k$ elements in the core-sparse matrix in~\Cref{tbl:mmf_cur_sparse}. 
\begin{table*}[!htb]
    \small 
    \centering 
    \begin{tabular}{lcccccc}
        \toprule 
        Matrix kind & num. matrices & 1\% & 10\% & 25\% & 50\% & 75\% \\
        \midrule
        \multicolumn{7}{c}{\textit{problems with underlying 2D/3D geometry}}\\
        structural problem&4&100.0&100.0&100.0&50.0&50.0\\
        materials problem&1&100.0&100.0&100.0&100.0&100.0\\
        2D/3D problem&5&100.0&100.0&100.0&80.0&60.0\\
        computational fluid dynamics problem&3&100.0&100.0&100.0&0.0&0.0\\
        \midrule
        \multicolumn{7}{c}{\textit{problems with no underlying geometry}}\\
        directed graph&18&100.0&50.0&22.3&16.7&0.0\\
        chemical process simulation problem&21&100.0&90.5&67.7&57.1&38.1\\
        circuit simulation problem&2&100.0&100.0&100.0&50.0&0.0\\
        power network problem&3&100.0&100.0&100.0&66.7&33.3\\
        counter-example problem&1&100.0&100.0&100.0&100.0&100.0\\
        economic problem&10&100.0&100.0&100.0&70.0&50.0\\
        statistical/mathematical problem&1&100.0&100.0&100.0&100.0&0.0\\
        \midrule
        total wins: & &100.0 &84.1 &69.6&49.3 &29.0 \\
        \bottomrule 
    \end{tabular}
    \caption{Matrix compression using asymmetric MMF (\textsc{GreedyTopN}). See~\Cref{tbl:mmf_cur_additive} for a description of the table entries. }
    \label{tbl:mmf_cur_sparse}
\end{table*}
Here, we see that MMF is able to capture structure and compress nearly all kinds of matrices at all compression levels.

From the numerical results, it seems that allowing interaction between wavelet levels is critical to achieving good approximation errors. This indicates that hierarchical structure in nonsymmetric matrices is more complex than that uncovered by simply extending the symmetric MMF analogous to an SVD. 

For smaller compression ratios, CUR seems to be preferable for some types of problems such as directed graphs. It is possible that some matrices may have global-scale behavior that is superimposed over the hierarchical structure. This is possible, for instance, when the data is noisy. In the next part, we explore ways of combining MMF and low-rank approaches to achieve even better compression errors. 

\subsubsection{Matrix compression via a combination of CUR and MMF}\label{sec:combination}

\begin{figure}[htb]
    \includegraphics[width=\linewidth]{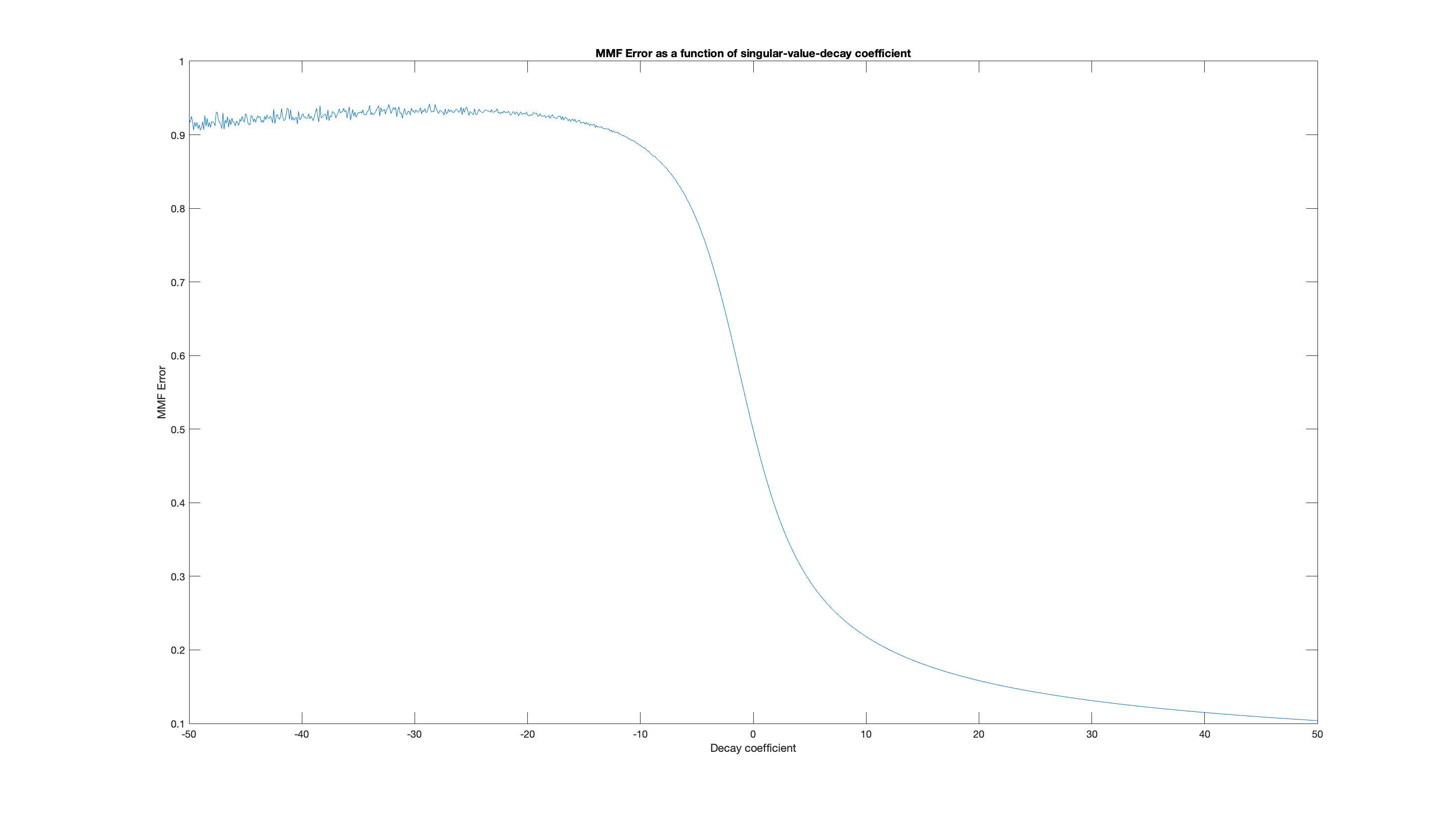}
    \caption{Approximation error of symmetric MMF with different rates of decay in the eigenvalues. The $200 \times 200$ random matrices were generated by computing $Q D Q^T$, where $Q\in\R^{200\times 200}$ is a fixed, randomly generated orthogonal matrix. $D$ is a diagonal matrix whose entries are given by $\frac{1-e^{t(x-1)}}{1-e^{t}}$, where $t$ is the decay coefficient and $x$ takes $200$ uniformly spaced values between $0$ and $1$. The X-axis corresponds to $t$.}
    \label{fig:mmferrors}
\end{figure}

In order to better motivate this section, we first consider~\Cref{fig:mmferrors}, which shows the error of MMF approximation of a matrix $A$ as a function of the decay rate of its singular values. 

\begin{figure}[htb]
    \includegraphics[width=\linewidth]{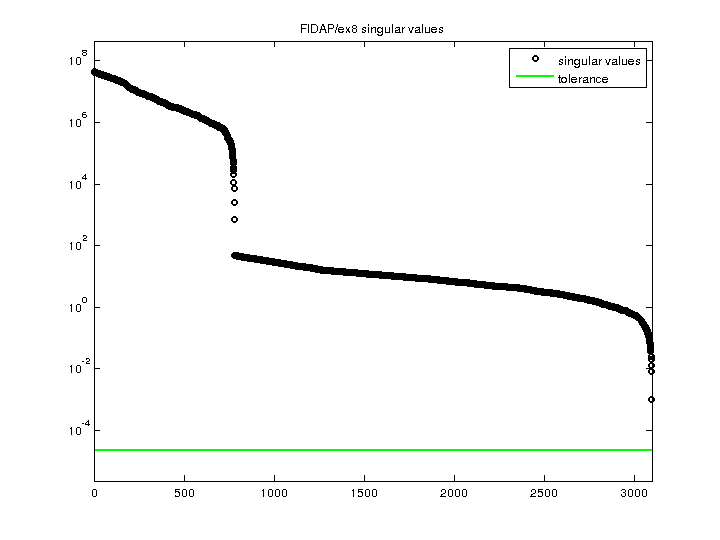}
    \caption{Plot showing the singular values of the `ex8' matrix from the UFlorida sparse matrix repository in descending order. Although the matrix is full-rank, the sharp drop in the singular values around index 800 imparts a low-rank nature to the matrix. (Note that the y-axis is on log-scale). Image sourced from UFlorida Sparse Matrix Repository~\citep{davis2011university}.}
    \label{fig:ex18}
\end{figure}

We observe that when the singular values decay at a super-linear rate, MMF does not approximate the matrix well. However, in many real-world matrices, the decay of singular values may not be uniform and may fluctuate to very high and low values. An example of such a spectrum is shown in~\Cref{fig:ex18}. In this case we see that although there are segments where the decay rate is sub-linear, the large decay rate around singular value 1500 may significantly impact the effectiveness of MMF on this matrix. Thus, it might be well motivated to first compress the matrix using low-rank methods and then use MMF on the residual. We describe this technique in~\Cref{alg:mmf_cur_combined}.

\begin{algorithm}
    \caption{Matrix compression using a combination of low-rank and multiresolution methods}
    \label{alg:mmf_cur_combined}
    \begin{algorithmic}
        \REQUIRE Square matrix $A$, compression size $k$, and $r$
        \STATE [C,U,R] = CUR(A, r)
        \STATE [P, H, Q] = MMF(CUR, k)
        \RETURN $\norm{PHQ^T}_F/\norm{A}_F$ 
    \end{algorithmic}
\end{algorithm}

\begin{figure}[htb]
    \centering 
    \includegraphics[width=0.45\linewidth]{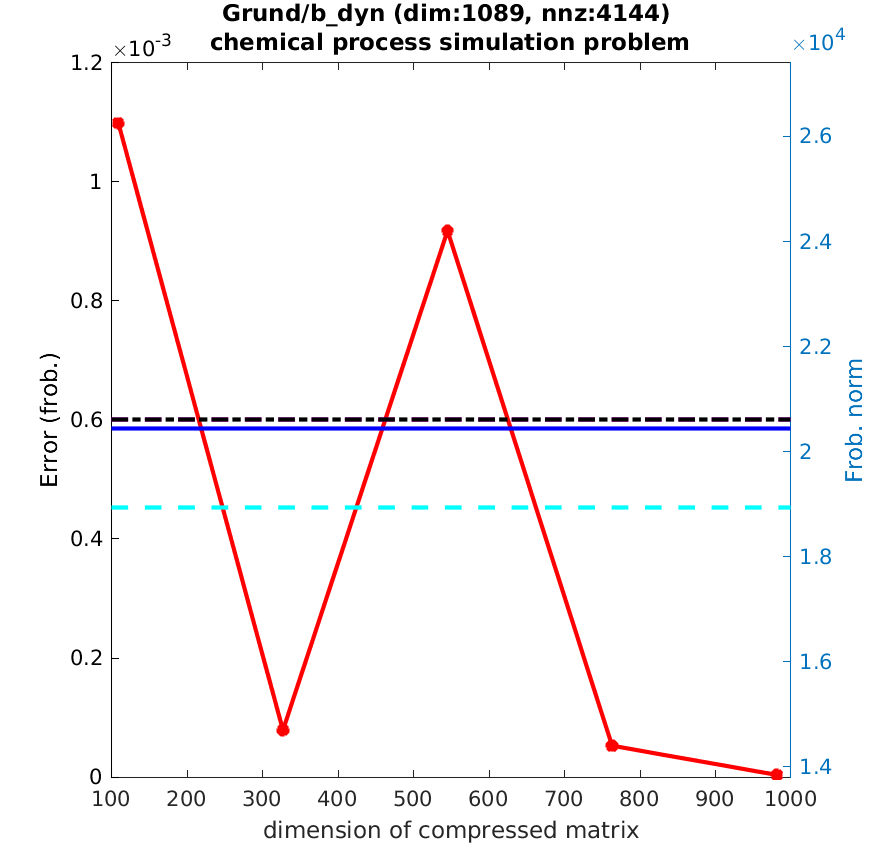} 
    \includegraphics[width=0.45\linewidth]{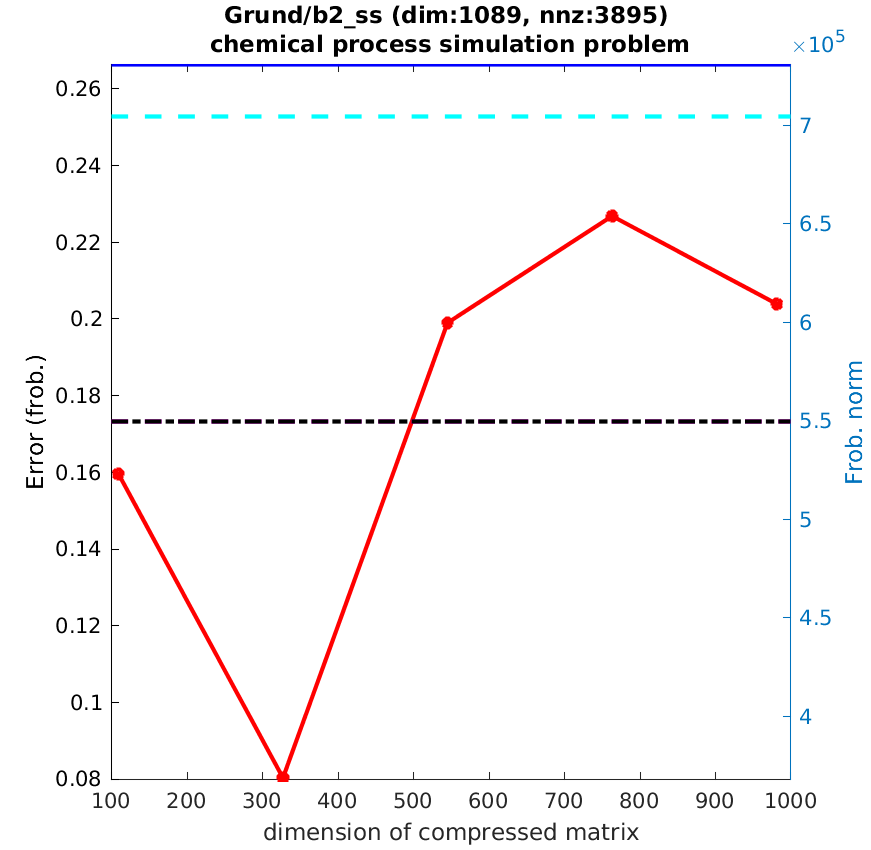} \\
    \includegraphics[width=0.45\linewidth]{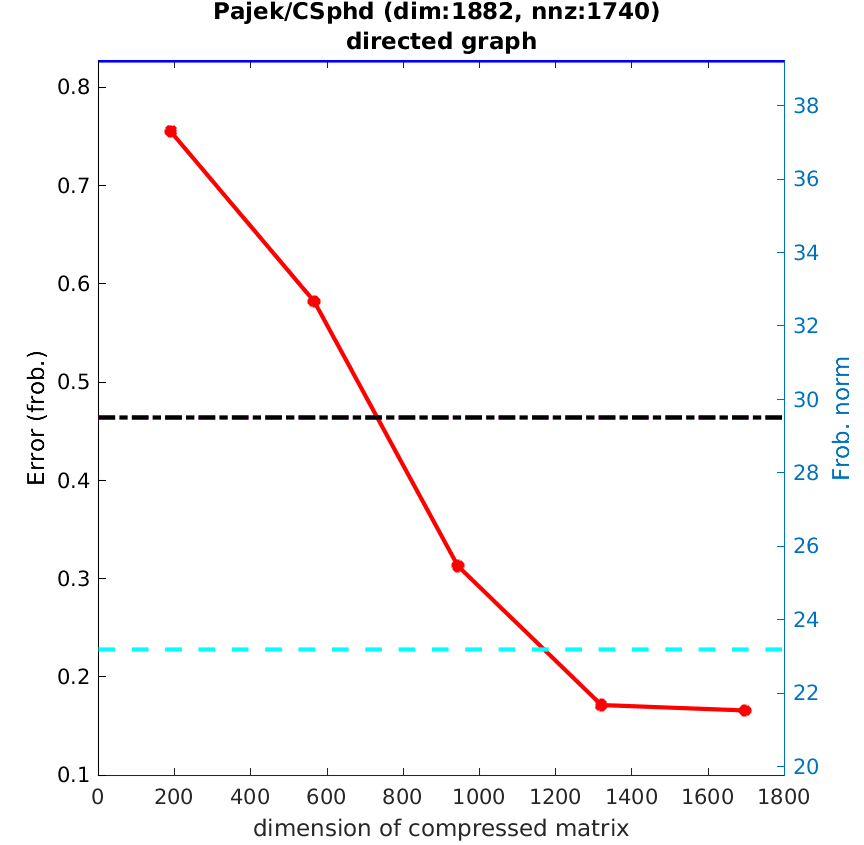}  
    \includegraphics[width=0.45\linewidth]{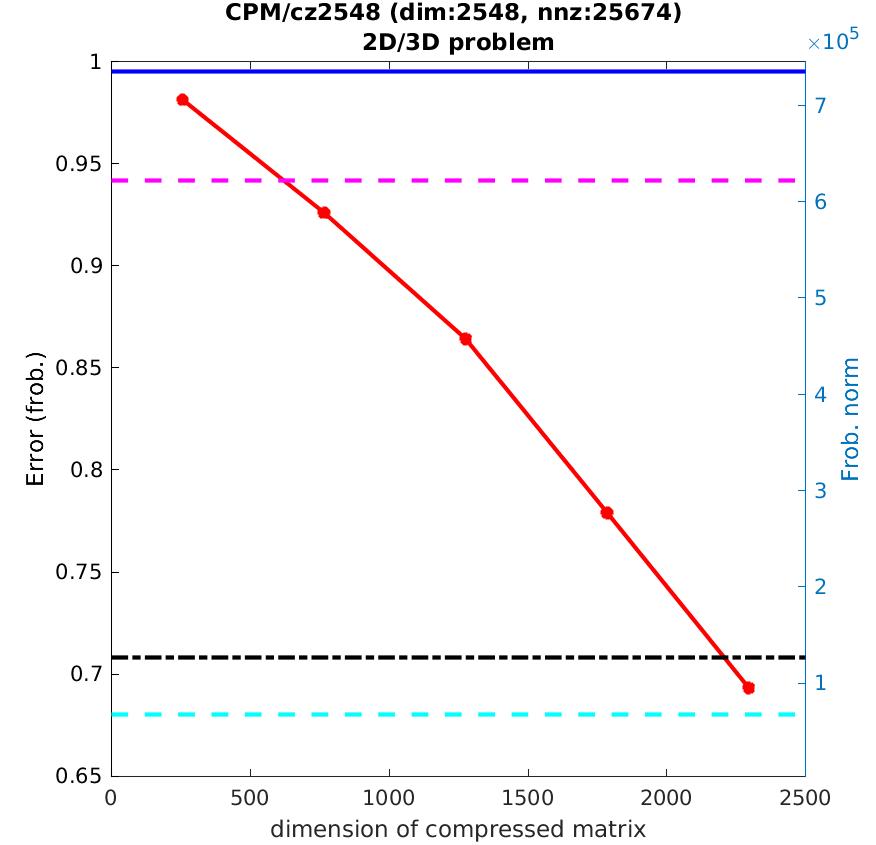} \\
    \includegraphics[width=0.45\linewidth]{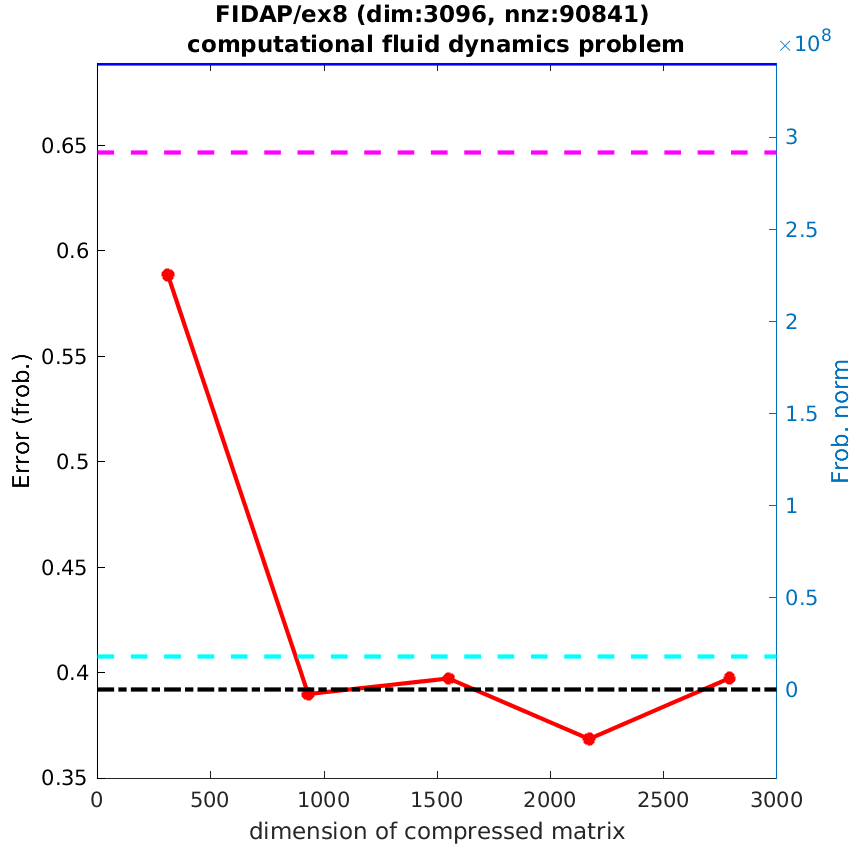} 
    \includegraphics[width=0.45\linewidth]{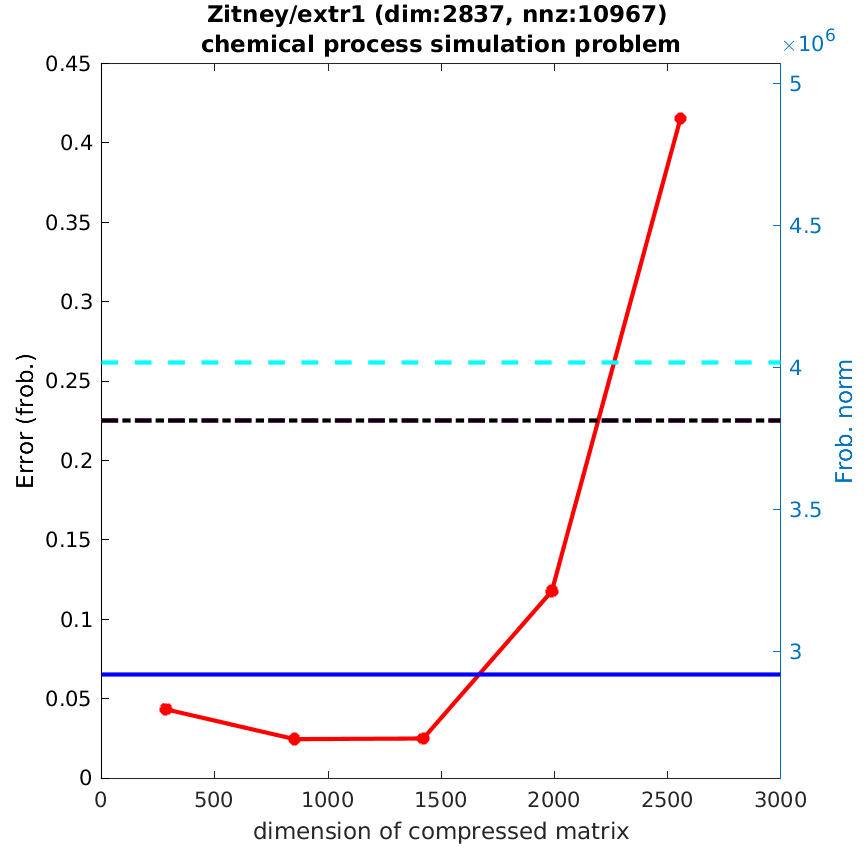} \\
    \caption{Frobenius norm error as a function of CUR compression rank for compressing the matrix to 5\% of its original size using the combined low-rank + multiresolution method. Solid blue line indicates error of compression with only CUR to 5\% matrix size, while the dotted blue line indicates the same, but with MMF.}
    \label{fig:mmf_cur_combined}
\end{figure}

For various values of $k$, we report the errors for compressing various matrices to 5\% of their original size in~\Cref{fig:mmf_cur_combined}.

In most of the matrices, we observe that there exists at least one rank to which the matrix can be first compressed using CUR before compressing with MMF, that results in lower compression errors than each of CUR and MMF methods individually. In practice, there is no method to determine what this rank would be for a given matrix, however, it is an interesting and worthwhile problem to be able to remove a small global-scale component of the matrix, after which the hierarchical structure may be easily extracted.

In all the experiments, the running time of MMF is equal or slightly slower than that of low-rank-based methods. However, we note that our algorithms can be made significantly faster by straightforwardly adapting parallel versions of MMF~\citep{kondor2015parallel} to the asymmetric case. We expect speedups similar to those achieved in the symmetric case~\citep{teneva2016multiresolution}.
\section{Conclusion}
In this paper, we have presented approaches to extend Multiresolution Matrix Factorization to nonsymmetric matrices. In general, efficiently estimating hierarchical structure in matrices is an increasingly important problem for which there currently exist very few practically useful methods. For future work we aim to use ideas from wavelet theory on more general spaces to define MMFs for general square matrices in which the MRA interpretation is clear.

%\textbf{Future work}: Note that \[\sup_{v\in V\subseteq\R^n} \frac{\inner{v,Av}}{\inner{v,v}} = \sup_{v\in V\subseteq\R^n} \frac{\inner{v,(A+A^T)v}}{2\inner{v,v}}\]
%Thus, performing a multiresolution analysis (\Cref{fig:space_splitting}) on $\R^n$ with respect to $A$ is identical as with $\frac{A+A^T}{2}$ which is a symmetric matrix. Suppose $U\in\R^{n\times n}$ represents the wavelet matrix from this multiresolution analysis. When $A$ is symmetric,~\cite{kondor2014multiresolution} defined the structure of $U^TAU$ as core-diagonal. In future research, it may be worthwhile to understand what kind of structure makes sense for the asymmetric case.
\bibliography{references}

\end{document}